\documentclass[12pt]{article}
\usepackage[a4paper, total={7in, 10in}]{geometry}
\usepackage[parfill]{parskip}
\usepackage{physics, tensor, float, subcaption}
\usepackage{graphicx}
\graphicspath{ {Plots/} }
\usepackage{jhep-mod}
\usepackage{bm}
\usepackage{soul}
\usepackage{amssymb,amsmath,amsthm}
\usepackage{mathrsfs}
\usepackage[utf8]{inputenc}
\usepackage{enumerate}
\usepackage{bigints}
\usepackage{xcolor}
\usepackage{appendix}
\usepackage{graphicx}
\usepackage{float}
\usepackage{tikz}
\usepackage{setspace}
\usepackage{cancel}
\usepackage{doi}
\definecolor{purple}{rgb}{1,0,1}
\definecolor{lime}{HTML}{A6CE39} % needs xcolor
\definecolor{lime}{HTML}{A6CE39}
\newcommand{\orcidicon}{%
	\begin{tikzpicture}
	\draw[lime, fill=lime] (0,0) 
		circle [radius=0.16] 
		node[white] {{\fontfamily{qag}\selectfont \tiny ID}};
	\draw[white, fill=white] (-0.0625,0.095) 
		circle [radius=0.007];
	\end{tikzpicture}
	\hspace{-5mm}
}
\newcommand\orcidMatt{{\href{https://orcid.org/0000-0003-1088-6485}{\orcidicon}}}
\renewcommand{\O}{\mathcal{O}}
\begin{document}
%========================================================
%========================================================
%========================================================

\title{
\leftline{Effective de la Vall\'e Poussin style bounds}
on the first Chebyshev function\\
}

%========================================================
%========================================================
%========================================================
\author{
\Large
Matt Visser\!\orcidMatt\!
}
%========================================================
%========================================================
%========================================================
%========================================================
\affiliation{School of Mathematics and Statistics, Victoria University of Wellington, \\
\null\qquad PO Box 600, Wellington 6140, New Zealand.}
%========================================================
%========================================================
\emailAdd{matt.visser@sms.vuw.ac.nz}
%========================================================
%========================================================
\def\theta{\vartheta}
\def\O{{\mathcal{O}}}

\abstract{
\vspace{1em}

In 1898 Charles Jean de la Vall\'e Poussin, as part of his famed proof of the prime number theorem, developed an ineffective bound on the first Chebyshev function of the form:\\
\[
 |\theta(x)-x| = \O\left(x \exp(-K \sqrt{\ln x})\right).
\]
This bound holds for $x$ sufficiently large, $x\geq x_0$, and $K$ some unspecified positive constant. 
To the best of my knowledge this bound has never been made effective --- I have never yet seen this bound made fully explicit, with precise values being given for $x_0$ and $K$. 
Herein, using a number of effective results established over the past  50 years, I shall develop two very simple explicit fully effective bounds of this type:
\[
|\theta(x)-x| <  \;
{x} \;\exp\left( - {1\over4} \sqrt{\ln x}\right);
\qquad (x\geq 2).
\]
\[
|\theta(x)-x| <  \;
{x} \;\exp\left( - {1\over3} \sqrt{\ln x}\right);
\qquad (x\geq 3).
\]
Many other fully explicit bounds along these lines can easily be developed. \\
For instance one can trade off stringency against range of validity: 
\[
|\theta(x)-x| <  \;
{1\over 2} \; {x} \;\exp\left( - {1\over4} \sqrt{\ln x}\right);
\qquad (x\geq 29),
\]
\[
|\theta(x)-x| <  \;
{1\over 2} \; {x} \;\exp\left( - {1\over3} \sqrt{\ln x}\right);
\qquad (x\geq 41).
\]
With hindsight, some of these effective bounds could have been established almost 50 years ago.

\bigskip

\bigskip
\noindent
{\sc Date:} 2 November 2022; \LaTeX-ed \today


}

%========================================================
\maketitle
%========================================================
\def\tr{{\mathrm{tr}}}
\def\diag{{\mathrm{diag}}}
\def\cof{{\mathrm{cof}}}
\def\pdet{{\mathrm{pdet}}}
\def\d{{\mathrm{d}}}
\parindent0pt
\parskip7pt
\def\Kerr{{\scriptscriptstyle{\mathrm{Kerr}}}}
\def\eos{{\scriptscriptstyle{\mathrm{eos}}}}
%================================================
\section{Introduction}
%================================================

In 1898 Charles Jean de la Vall\'e Poussin developed an ineffective bound on the first Chebyshev function of the form~\cite{Poussin}:
\begin{equation}
\label{E:Poussin}
|\theta(x)-x| = \O\left(x \exp(-K \sqrt{\ln x})\right).
\end{equation}
This bound holds for $x$ sufficiently large, $x\geq x_0$, and $K$ some unspecified positive constant. 

Subsequent work over the last 50 years has developed a large number of related but distinct fully effective bounds of the form~\cite{Schoenfeld,Trudgian,Johnston-Yang,Fiori-et-al}:
\begin{equation}
|\theta(x)-x| <  a \;x \;(\ln x)^b \; \exp\left(-c\; \sqrt{\ln x}\right); 
\qquad (x \geq x_0).
\end{equation}
\begin{itemize}
\item 
For some widely applicable effective bounds of this type see Table I.\\
(A straightforward elementary numerical computation is required to determine the numerical coefficients in the Schoenfeld~\cite{Schoenfeld} and Trudgian~\cite{Trudgian} bounds.)
\item
For some asymptotically more stringent effective bounds of this type, but valid on significantly more restricted regions, see Table~II (based on reference~\cite{Johnston-Yang}), and the extensive tabulations in reference~\cite{Broadbent-et-al}.
\end{itemize}
\enlargethispage{30pt}
What I have not yet seen is any attempt to take the effective bounds of Tables I and II and use them to 
make the original de la Vall\'e Poussin bound fully effective. 
Here are two particularly clean fully effective versions of the de la Vall\'e Poussin bound:
\begin{equation}
|\theta(x)-x| <  \;
{x} \;\exp\left( - {1\over4} \sqrt{\ln x}\right);
\qquad (x\geq 2).
\end{equation}
\begin{equation}
|\theta(x)-x| <  \;
{x} \;\exp\left( - {1\over3} \sqrt{\ln x}\right);
\qquad (x\geq 3).
\end{equation}
I shall explain how to derive these bounds below. 
\clearpage

%=================================================
\begin{table}[!h]
\caption{Some widely applicable effective bounds.}\smallskip
\begin{center}
\begin{tabular}{||c|c|c|c||c||}
\hline
\hline
$a$ & $b$ & $c$ & $x_0$ & Source \\
\hline
\hline
 0.2196138920& 1/4 & 0.3219796502 & 101 & Schoenfeld~\cite{Schoenfeld}\\
\hline
\hline
0.2428127763 & 1/4 &0.3935970880 & 149 & Trudgian~\cite{Trudgian}\\
\hline
\hline
9.220226 & 3/2 & 0.8476836 & 2 & Fiori--Kadiri--Swidinsky~\cite{Fiori-et-al}\\
\hline
\hline
9.40 & 1.515 & 0.8274 & 2 & Johnston--Yang~\cite{Johnston-Yang}\\
\hline
\hline 
\end{tabular}
\end{center}
\end{table}
%=================================================
 
%=============================================
\begin{table}[!htb]
\caption{Asymptotically stringent effective bounds valid on restricted regions~\cite{Johnston-Yang}.}\smallskip
\begin{center}
\begin{tabular}{||c|c|c|c||}
\hline
\hline
$a$ & $b$ & $c$ & $x_0$  \\
\hline
\hline
8.87 & 1.514 & 0.8288 & $\exp(3000)$ \\
8.16 & 1.512 & 0.8309 & $\exp(4000)$  \\
7.66 & 1.511 & 0.8324 & $\exp(5000)$  \\
7.23 & 1.510 & 0.8335 & $\exp(6000)$  \\
7.00 & 1.510 & 0.8345 & $\exp(7000)$  \\
6.79 & 1.509 & 0.8353 & $\exp(8000)$  \\
6.59 & 1.509 & 0.8359 & $\exp(9000)$  \\
6.73 & 1.509 & 0.8359 & $\exp(10000)$ \\
\hline\hline
23.14 & 1.503 & 0.8659 & $\exp(10^5)$  \\
38.58 & 1.502 & 1.0318 & $\exp(10^6)$  \\
42.91 & 1.501 & 1.0706 & $\exp(10^7)$  \\
44.42 & 1.501 & 1.0839 & $\exp(10^8)$  \\
44.98 & 1.501 & 1.0886 & $\exp(10^9)$  \\
45.18 & 1.501 & 1.0903 & $\exp(10^{10})$  \\
\hline
\hline
\end{tabular}
\end{center}
\end{table}
%=================================================

%================================================
\section{Strategy}
%================================================
%================================================

Note that for any $b>0$, $c>0$, and any $\tilde c\in(0,c)$, elementary calculus implies:
\begin{eqnarray}
(\ln x)^b \; \exp\left(-c\; \sqrt{\ln x}\right) 
&=& \left\{ (\ln x)^b \exp\left(-[c-\tilde c] \; \sqrt{\ln x}\right) \right\} \exp\left(-\tilde c\; \sqrt{\ln x}\right) 
\nonumber\\
& \leq & \left\{ \left( 2b\over c-\tilde c \right)^{2b} \exp(-2 b)  \right\} \exp\left(-\tilde c\; \sqrt{\ln x}\right).
\end{eqnarray}
The key observation here is that the quantity in braces is explicitly bounded, 
and achieves a global maximum at 
$x_{peak} = \exp\left( \left[2b/(c-\tilde c)\right]^2\right)$.
Consequently we have the following lemma.

\clearpage

\paragraph{Lemma:}
\emph{Any effective bound of the form
\begin{equation}
|\theta(x)-x| <  a \;x \;(\ln x)^b \; \exp\left(-c\; \sqrt{\ln x}\right); 
\qquad (x \geq x_0).
\end{equation}
implies the existence of another effective bound of the de la Vall\'e Poussin form
\begin{equation}
|\theta(x)-x| <  \tilde a \; x\; \exp\left(-\tilde c\; \sqrt{\ln x}\right); 
\qquad (x \geq x_*; \;\; x_*  \leq x_0).
\end{equation}
Here  $\tilde c$ is an arbitrary number in the interval $\tilde c \in (0,c)$ and 
\begin{equation}
\tilde a = a \left( 2 b\over c-\tilde c \right)^{2b} \exp(-2 b).
\end{equation}
Note that this new bound of the de la Vall\'e Poussin form certainly holds for $x>x_0$, but if $x_0$ is sufficiently small one might be able to widen the range of applicability to some new $x \geq x_*$ with $x_* \leq x_0$ by explicit computation. 
}

We now apply this lemma to the various bounds explicated above. 

%================================================
\section{Some effective bounds}
%================================================

\enlargethispage{20pt}
First let us consider some widely applicable derived bounds of the de la Vall\'e Poussin form, as presented in Table~III. 
Note that the selection of a specific value of $\tilde c$ is a \emph{choice}, and the computation of $\tilde a$ is then immediate --- there is an infinite number of other 
effective 
bounds of de la Vall\'e Poussin form that we could develop.
Determining $x_*$ then requires computationally checking low values of $x$. 

%=====================================================
\begin{table}[!htb]
\caption{Some widely applicable derived bounds of the de la Vall\'e Poussin form.}\smallskip
\begin{center}
\begin{tabular}{||c|c|c||c||}
\hline
\hline
$\tilde a$ & $\tilde c$ & $x_*$ & Based on \\
\hline
\hline
 0.3510691792& 1/4  & 59 & Schoenfeld~\cite{Schoenfeld}\\
\hline
\hline
0.2748124978 & 1/4 & 101 & Trudgian~\cite{Trudgian}\\
\hline
0.4242102935 & 1/3 & 59 & Trudgian~\cite{Trudgian}\\
\hline
\hline
295 & 1/2 & 2 & Fiori--Kadiri--Swidinsky~\cite{Fiori-et-al}\\
\hline
\hline
385 & 1/2  & 2 & Johnston--Yang~\cite{Johnston-Yang}\\
\hline
\hline
\end{tabular}
\end{center}
\end{table}
%===============================================

By now relaxing the prefactor $\tilde a$, one can increase the range of validity of the bound, (ie, decrease $x_*$). In this way, after some computation, one finds
\begin{equation}
|\theta(x)-x| <  \;
{x} \;\exp\left( - {1\over4} \sqrt{\ln x}\right);
\qquad (x\geq 2).
\end{equation}
Note that this could in principle have been deduced as early as 1976, some 46 years ago,
from the work of Schoenfeld~\cite{Schoenfeld}.

Similarly
\begin{equation}
|\theta(x)-x| <  \;
{x} \;\exp\left( - {1\over3} \sqrt{\ln x}\right);
\qquad (x\geq 3).
\end{equation}
Note that this particular bound could in principle have been deduced as early as  2016, 
some 6 years ago, from the work of Trudgian~\cite{Trudgian}. 
For the other two widely applicable bounds, some preliminary experimental investigations \emph{suggest} that it might be possible to reduce the numerical prefactors (295, 385) significantly --- but doing so would require rather different techniques from the elementary observations made above. 

As always one can trade of stringency against range of validity. In this regard let me mention two specific examples 
\begin{equation}
|\theta(x)-x| <  \;
{1\over 2} \; {x} \;\exp\left( - {1\over4} \sqrt{\ln x}\right);
\qquad (x\geq 29),
\end{equation}
and
\begin{equation}
|\theta(x)-x| <  \;
{1\over 2} \; {x} \;\exp\left( - {1\over3} \sqrt{\ln x}\right);
\qquad (x\geq 41).
\end{equation}

In contrast, for some asymptotically stringent bounds, based on the Johnston--Yang results presented in reference~\cite{Johnston-Yang}, consider Table~IV. Note that we can again make the exponential factor smaller, (by increasing $\tilde c$), at the cost of making the numerical prefactor $\tilde a$ larger. 
For instance one can deduce the ineffective bound
\begin{equation}
|\theta(x)-x| =  \;
\O \left( {x} \;\exp\left( -  \sqrt{\ln x}\right)\right),
\end{equation}
which can be made effective as (for instance):
\begin{equation}
|\theta(x)-x| =  \;
83063 \;\; {x} \;\exp\left( -  \sqrt{\ln x}\right); \qquad \left(x> \exp(10^{10})\right).
\end{equation}
Many variations on this theme can be developed. 

%=====================================================
\begin{table}[!htb]
\caption{Asymptotically stringent bounds of the de la Vall\'e Poussin form.}\smallskip
\begin{center}
\begin{tabular}{||c|c|c||}
\hline
\hline
$\tilde a$ & $\tilde c$ & $x_*$  \\
\hline
\hline
357 & 1/2  & $\exp(3000)$\\
320 & 1/2  & $\exp(4000)$  \\
295 & 1/2  & $\exp(5000)$  \\
274 & 1/2  & $\exp(6000)$  \\
263 & 1/2 &  $\exp(7000)$  \\
252 & 1/2  & $\exp(8000)$  \\
244 & 1/2 & $\exp(9000)$  \\
249 & 1/2  & $\exp(10000)$  \\
\hline
\hline
644 & 1/2  & $\exp(10^5)$  \\
348 & 1/2 & $\exp(10^6)$  \\
312 & 1/2 & $\exp(10^7)$  \\
301 & 1/2 &  $\exp(10^8)$  \\
298 & 1/2 &  $\exp(10^9)$  \\
297 & 1/2 & $\exp(10^{10})$  \\
\hline
\hline
1642333 & 1 & $\exp(10^6)$  \\
165152 & 1 & $\exp(10^7)$  \\
101831 & 1 &  $\exp(10^8)$  \\
87551 & 1 &  $\exp(10^9)$  \\
83063 & 1 & $\exp(10^{10})$  \\
\hline
\hline
\end{tabular}
\end{center}
\label{default}
\end{table}%

%================================================
%================================================
%================================================
\section{Conclusions}\label{S:discussion}
%================================================

With some hindsight, deriving effective bounds of the de la Vall\'e Poussin form is, (given various effective results obtained over the last 50 years~\cite{Schoenfeld,Trudgian,Johnston-Yang,Fiori-et-al,Broadbent-et-al}), seen to be almost trivial.  
Certainly, (given the effective results reported in~\cite{Schoenfeld,Trudgian,Johnston-Yang,Fiori-et-al,Broadbent-et-al}), nothing deeper than elementary calculus and some slightly tedious numerical checking was required. On the other hand, conceptually it is very pleasant to see simple explicit and effective bounds of the de la Vall\'e Poussin form dropping out so nicely.

%================================================
\section*{Acknowledgements}
%================================================

MV was directly supported by the Marsden Fund, \emph{via} a grant administered by the
Royal Society of New Zealand.

%================================================
%\appendix
%================================================
%\section{}
%=====================================================
%\label{A:A}

\clearpage
%================================================
%================================================

%================================================
\end{document}